\documentclass[12pt]{amsart}
\usepackage{amssymb, mathrsfs,amsmath}
\usepackage{latexsym}
\usepackage{amsfonts}
\usepackage{shadow}
\newcommand{\longhookrightarrow}
{\ensuremath{\lhook\joinrel\relbar\joinrel\relbar\joinrel\rightarrow}}
\newcommand{\w}{\omega}
\date{}
\author[D. Alpay]{Daniel Alpay}
\author[P. Jorgensen]{Palle Jorgensen}
\address{(DA) Department of Mathematics \newline
Ben Gurion University of the Negev \newline P.O.B. 653, \newline
Be'er Sheva 84105, \newline ISRAEL} \email{dany@math.bgu.ac.il}
\address{(PJ)
Department of Mathematics\newline 14 MLH \newline The University
of Iowa Iowa City,\newline IA 52242-1419 USA}
\email{jorgen@math.uiowa.edu}
\thanks{D. Alpay thanks the
Earl Katz family for endowing the chair which supported his
research. This research was supported in part by the Binational
Science Fundation grant number 2010117}

\title[Stochastic processes induced by singular operators]
{Stochastic processes induced by singular operators}
\begin{document}
\parindent 0cm
\newtheorem{Pa}{Paper}[section]
\newtheorem{Tm}[Pa]{{\bf Theorem}}
\newtheorem{La}[Pa]{{\bf Lemma}}
\newtheorem{Cy}[Pa]{{\bf Corollary}}
\newtheorem{Rk}[Pa]{{\bf Remark}}
\newtheorem{Pn}[Pa]{{\bf Proposition}}
\newtheorem{Dn}[Pa]{{\bf Definition}}
\newtheorem{Ex}[Pa]{{\bf Example}}
\numberwithin{equation}{section}
\def\L{\mathbf L}
\def\R{\mathbb R}
\def\N{\mathbb N}
\def\C{\mathbb C}
\def\s{\mathscr S}
\def\H{\mathscr H}
\def\ss{\mathscr S^\prime}
\def\sr{\mathscr S(\R)}
\def\ssr{\mathscr S'(\R)}
\def\(s){\mathscr S(\R^2)}
\def\F{\mathcal F}
\def\P{\mathcal P}
\def\W{\mathcal W}
\def\Dom{{\rm dom}~(T_m)}
\def\Doms{{\rm dom}~(T_m^*)}
\def\Def{\stackrel{{\rm def.}}{=}}
\begin{abstract}
In this paper we study a general family of multivariable Gaussian
stochastic processes. Each process is prescribed by a fixed Borel
measure $\sigma$ on $\mathbb R^n$. The case when $\sigma$ is
assumed absolutely continuous with respect to Lebesgue measure
was studied earlier in the literature, when $n=1$. Our focus here
is on showing how different equivalence classes (defined from
relative absolute continuity for pairs of measures) translate
into concrete spectral decompositions of the corresponding
stochastic processes under study. The measures $\sigma$ we
consider are typically purely singular.
Our proofs rely on the theory of (singular) unbounded operators in
Hilbert space, and their spectral theory.
\end{abstract}
 \maketitle \tableofcontents
\section{Introduction}
\setcounter{equation}{0}
We study Gaussian stochastic processes with stationary
mean-square-increments. These are used in for example modeling of
time-dependent phenomena in signal processing, in information
theory, in telecommunication, and in a host of applications. A
process $\left\{X(t)\right\}$ indexed by $t\in\mathbb R^n$ on a
probability space is said to have stationary mean-square
increments if the mean-square of each of the increments
$X(t)-X(s)$ depends only on the difference $t-s$. Note that we are
not restricting the nature of the statistical distribution of the
increments, only the expectation of the mean-square. In
particular, we are not assuming independent increments.
Nonetheless, this mean-square property allows us to analyze
decompositions of the processes and stochastic integrals with the
use of spectral theory and tools from the theory of operators in
Hilbert space.\\

There are several reasons for the importance of stationary
mean-square increment processes. First they include a rich family
of processes from applications; secondly, these are precisely the
processes that admit a canonical  action  of a unitary group of
time-transformations; the unitary operators act on
$\mathbf L_2(\Omega, \mathcal A,
P)$ where $(\Omega, \mathcal A,P)$ is the underlying
probability space. For $n=1$, it will be a one-parameter group of
unitary operators; and in general, $n>1$, it will be a unitary
representation of the additive group $\mathbb R^n$. In the case
of $n=1$, and Brownian motion, Wiener and Kakutani proved that
this one-parameter group is a group of point transformations
acting ergodically; see \cite{Kak48}.
The third reason is purely mathematical: the
class of stationary mean-square increments processes admits a
classification. We turn to this in Theorem \ref{thm:central}
below. In Section \ref{sec:ergodicity} we introduce the unitary
representations of $\mathbb R^n$.\\

It is possible to characterize the processes mentioned in the
previous paragraphs by the family $\mathcal C$ of regular positive
Borel measures $\sigma$ on $\mathbb R^n$ subject to
\begin{equation}
\label{eq:bound} \int_{\mathbb
R^n}\frac{d\sigma(u)}{(1+|u|^2)^p}<\infty
\end{equation}
for some $p\in\mathbb N_0$.
Such a measure $\sigma$ is the spectral
function of a homogeneous generalized stochastic field in the
sense of Gelfand. See \cite[p. 283]{MR35:7123}. Here we study
this correspondence in reverse. Earlier work  so far was
restricted to the case when the spectral density is assumed to be
absolutely continuous and $n=1$, see \cite{aal2}. The case of
singular measure was considered in \cite{MR2793121}. Here, we are
making no restriction on the spectral type.
We associate to $\sigma\in\mathcal C$ four natural objects:\\

$(a)$ The quadratic form
\begin{equation}
\label{qsigma}
q_\sigma(\psi)=\int_{\mathbb
R^n}|\widehat{\psi}(u)|^2d\sigma(u),
\end{equation}
where $\psi$ belongs to the Schwartz space $\mathcal S_{\mathbb
R^n}$, and where
\[
\widehat{\psi}(u)=\int_{\mathbb R}e^{-iux}\psi(x)dx
\]
denotes the Fourier transform of $\psi$.\\

$(b)$ A linear operator $Q_\sigma$ such that
\[
q_\sigma(\psi)=\|Q_\sigma\psi\|^2_{{\mathbf L_2(\mathbb
R^n,dx)}},\quad \psi\in\mathcal S_{\mathbb R^n}.
\]

$(c)$ A generalized stochastic process $\left\{X_\sigma(\psi),
\psi\in\mathcal S_{\mathbb R^n}\right\}$ such that
\[
E[X_\sigma(\psi_1)\overline{X_\sigma(\psi_2)}]=\int_{\mathbb R^n}
\widehat{\psi_1}(u)\overline{\widehat{\psi_2}(u)}d\sigma(u),
\]
for $\psi_1,\psi_2\in\mathcal S_{\mathbb R^n}$.\\

$(d)$ Let $\mathcal S_{\mathbb R^n}(\mathbb R)$ the Schwartz
space of {\sl real-valued} Schwartz functions of $n$-variables,
and
\begin{equation}
\label{eq:Omega} \Omega=\mathcal S^\prime_{\mathbb R^n}(\mathbb R)
\end{equation}
its dual, the space of all tempered distributions. The fourth object
associated to $\sigma$ is a
probability measure $d\mu_\sigma$ on $\mathcal S^\prime_{\mathbb
R^n}(\mathbb R)$ such that
\[
e^{-\frac{\|\widehat{\psi}\|^2_{\mathbf
L_2(d\sigma)}}{2}}=\int_{\Omega} e^{i\langle \w,\psi\rangle}
d\mu_\sigma(\w),\quad\psi\in\mathcal S_{\mathbb R^n}(\mathbb R),
\]
where $\langle\, ,\,\rangle$ denotes the duality between
$\mathcal S_{\mathbb R^n}(\mathbb R)$ and $\mathcal S_{\mathbb
R^n}^\prime(\mathbb R)$. One purpose of this paper is to study
various connections between
these quantities. See in particular Theorem \ref{thm:central}.\\

We note that the measures $\sigma$ of interest in building stochastic
processes include those for which $\mathbf L_2(\sigma)$ has a Fourier basis,
i.e.,
an orthogonal and total family of complex exponentials, where the frequencies
are made up from some discrete subset $S$ in $\mathbb R^n$. When this happens,
we say that $(\sigma, S)$ is a spectral pair. They have been extensively
studied in \cite{DuJo09, JoSo09a, JoSo09b, DJP09, DHJ09}.
Aside from their application to the constructing highly fluctuating
stochastic processes, the spectral pair measures have additional uses:
The measures $\sigma$  that are part of a spectral pair have the
following features:\\
$\mathbf 1.$
They are different from those usually considered in the theory of
stochastic processes, and yet they produce explicit models with
easy rules for computation.\\
$\mathbf 2.$
They possess intriguing multiscale properties, and selfsimilarity.\\

We also note that one can extend, as in our previous paper
\cite{MR2793121}, the generalized process $X_\sigma$ to be
defined for $t\in\mathbb R^n$.
When $n=1$, this is done by replacing the
variable $\phi\in\mathcal S$ by the functions
$\xi_t(u)=\frac{e^{itu}-1}{u}$.\\

For $n=1$, the quadratic form \eqref{qsigma} appears in
particular in the book of Gelfand and Vilenkin; see
\cite[Th\'eor\`eme $1^\prime$, p. 258]{MR35:7123}. Still for
$n=1$, the papers \cite{aal2} and \cite{MR2793121} can be seen as
the beginning of such a study, for certain absolutely continuous
measures, and for certain singular measures respectively. For
other treatments of the theory of Gaussian processes, we refer to
the following books, \cite{ItMc65}, a classic, and  \cite{HiHi93}
of a more recent vintage.

\begin{Rk}
We note that in the paper both the spaces of real-valued Schwartz
functions $\mathcal S_{\mathbb R^n}(\mathbb R)$ and
complex-valued Schwartz functions $\mathcal S_{\mathbb R^n}$, and
their duals, come into play. The real-valued case appears mostly
when we make use of the Bochner-Minlos theorem. (See for example
equations \eqref{wns} and \eqref{eq:iso} below.) And the
corresponding complex spaces are used in computation of adjoint
operators, with reference to some given Hilbert space inner
product. When needed, it is simple to insert a complexification.
\end{Rk}

The outline of the paper is as follows. The paper consists of
eleven sections besides the introduction. Section 2 is a
preliminary section where we review Hida's construction of the
white noise space. Given a positive measure $\sigma$
on $\mathbb R^n$, satisfying
\eqref{eq:bound},  we study in Sections \ref{sec:3} and
\ref{sec:4} the
associated quadratic form $q$ (:=
$q_\sigma$) defined by \eqref{qsigma}.
In Theorem \ref{tm:closable} we prove that $q$ is
closable if and only if $\sigma$ is absolutely continuous with
respect to Lebesgue measure on $\mathbb R^n$. These results have implications
on the nature of our stochastic processes $X_\sigma$
(indexed by the Schwartz space $\mathcal S_{\mathbb R^n}$), which is built in
Section 5. See  Theorem \ref{thm:central} where we construct the
associated process $X_\sigma$, and its path-space measure
$\mu_\sigma$ on the Schwartz space $\mathcal S^\prime_{\mathbb
R^n}(\mathbb R)$ of real-valued tempered distributions (see
\eqref{eq:Omega}).
In Theorem \ref{tm:iso}, we realize
$\mathbf L_2(\mathcal S^\prime_{\mathbb R^n}(\mathbb R),
\mathcal B(\mathcal S^\prime_{\mathbb R^n}(\mathbb R)), \mu_\sigma)$ in an explicit
manner as a reproducing kernel Hilbert space. A second generalized
stochastic process associated to $\sigma$ is constructed in Section 7.
Ergodicity is studied in Section 8.
In section \ref{sigma_function}, we prove
that two processes constructed from measures $\sigma_1$  and  $\sigma_2$
from $\mathcal C$ are orthogonal in the white noise space $\mathcal W$  if and
only $\sigma_1$  and  $\sigma_2$  are mutually singular measures.
More generally, we introduce a Hilbert space $\mathscr H$
of sigma-functions such
that each measure $\sigma$ from $\mathcal C$ generates a closed
subspace $\mathscr H(\sigma)$ in $\mathscr H$, with pairs of
mutually singular measures generating orthogonal subspaces in
$\mathscr H$.
In Theorem \ref{tm:9.1}, we show that our stochastic processes
$X_\sigma$ are indexed via an infinite-dimensional Fourier
transform by equivalence classes of measures $\sigma$, with
equivalence meaning mutually absolute continuity. As a result we
note that the case of closable $q$ in Theorem \ref{tm:closable}
only represents a single equivalence class of processes.
Therefore our study of the possibilities for measures $\sigma$
not absolutely continuous with respect to Lebesgue measure,
introduces a multitude of new processes not studied earlier. An
example is worked out in Section \ref{sec:12}; see also
\cite{MR2793121} for related examples. A generalized
Karhunen-Lo\`eve expansion is studied in Section 10. Our result,
Theorem \ref{tm:7.2}, offers a direct integral decomposition of
the most general stochastic process $X_\sigma$. Its conclusion
yields a formula for separation of variables (i.e., time and
sample-point) in the sense of Karhunen-Lo\`eve; decomposing
$X_\sigma$ as a direct integral of standard i.i.d. $N(0, 1)$
random variables. We prove in Section 11 that the space $\mathcal
C$ of measures satisfying \eqref{eq:bound} for some $p$ is a
module. As mentioned above, the last section is devoted to an
example.
\section{Preliminaries}
\setcounter{equation}{0}
There are generally two approaches to measures on function spaces
(path-space). The first one builds the measure from prescribed
transition probabilities and uses a limit construction, called
the Kolmogorov consistency principle. The second approach is based
on the notion of generating function (computed from a probability
measure), essentially a Fourier transform, obtained from the
Bochner-Minlos theorem. This theorem insures the existence of a
measure $d\mu_W$ on $\Omega=\mathcal S^\prime_{\mathbb
R^n}(\mathbb R)$ (as defined by \eqref{eq:Omega}) such that
\begin{equation}
\label{wns}
e^{-\frac{\|\psi\|^2_{{\mathbf L_2(\mathbb
R^n,dx)}}}{2}}=\int_{\Omega} e^{i\langle
\w,\psi\rangle}d\mu_W(\w),\quad \psi\in\mathscr S_{\mathbb
R^n}(\mathbb R).
\end{equation}

The white noise space is defined to be
\begin{equation}
\mathcal W=\mathbf L_2(\Omega,\mathcal B,d\mu_W),
\end{equation}
where $\mathcal B$ denotes the sigma-algebra of Borel sets.
Furthermore, \eqref{wns} induces an isomorphism $f\mapsto
\widetilde{f}$ from the space of square summable real-valued
functions of $\mathbf L_2(\mathbb R^n,dx)$ into $\mathcal W$ via
the formula
\begin{equation}
\label{eq:iso}
\widetilde{\psi}(\w)=\langle \w,\psi\rangle,\quad
\w\in\Omega,
\end{equation}
first for $\psi\in\mathcal S_{\mathbb R^n}(\mathbb R)$ and then by
continuity for every real-valued $f\in\mathbf L_2(\mathbb
R^n,dx)$. Throughout our paper, we will be using the
Bochner-Minlos' theorem for a variety of positive definite
functions on $\mathcal S_{\mathbb R^n}(\mathbb R)$. Each will be
required to be continuous on $\mathcal S_{\mathbb R^n}(\mathbb
R)$ with respect to the Fr\'echet topology. The positive definite
functions we consider will have a form similar to that in the
expression on the left hand side in \eqref{wns}, but with a
different quadratic form occurring in the exponent.
For background references we recommend \cite{MR1244577},
\cite{MR2444857}, \cite{new_sde} and \cite{MR1387829}.\\

\section{A sesquilinear form}
\setcounter{equation}{0}
\label{sec:3}
%
In our construction of Gaussian measures on the
space of real tempered distributions in a Gelfand-triple via the
Bochner-Minlos theorem, we will be making use of families of
sesquilinear forms on the Schwartz space $\mathcal S_{\mathbb
R^n}(\mathbb R)$. The properties of the measures in turn depend
on the nature of the sesquilinear forms under consideration, and
we now turn to these below. We first recall that a map from the
space of complex-valued Schwartz functions $\mathcal S_{\mathbb
R^n}$ into $\mathbb R_+$ is called a positive quadratic form if
for every $\phi,\psi\in \mathcal S_{\mathbb R^n}$ and
$c\in\mathbb C$, it holds that:
\begin{align}
\frac{1}{2}\left(q(\phi+\psi)+q(\phi-\psi)\right)&=q(\phi)+q(\psi),\\
q(c\psi)&=|c|^2q(\psi),\\
q(\psi)&\ge 0.
\end{align}

Let $u=(u_1,u_2,\ldots, u_n)$ and $v=(v_1,v_2,\ldots, v_n)$ be in $
\mathbb R^n$. We denote by
\[
u.v=\sum_{k=1}^n u_kv_k\quad{\rm and}\quad |u|^2=\sum_{k=1}^nu_k^2
\]
the inner product of $u$ and $v$ and the norm of $u$ respectively.

\begin{Dn}
\label{dn:C}
We denote by $\mathcal C$ the space of positive
measures $d\sigma$ on $\mathbb R^n$ such that \eqref{eq:bound}
holds for some $p\in\mathbb N$:
\begin{equation*}
\int_{\mathbb R^n}\frac{d\sigma(u)}{(1+|u|^2)^p}<\infty
\end{equation*}
\end{Dn}

We associate to $d\sigma$ the space
$\mathcal M(\sigma)$ of measurable functions $\psi$ such that
\[
\int_{\mathbb R^n}|\widehat{\psi}(u)|^2d\sigma(u)<\infty.
\]
In view of \eqref{eq:bound}, the space $\mathcal M(\sigma)$
contains the Schwartz space $\mathcal S_{\mathbb R^n}$ of
Schwartz functions of $n$ real variables. The set of functions
$\psi$ such that $q_\sigma(\psi)<\infty$ plays an important role
in \cite[Theorem 3.1, p. 258]{MR1790083} for $d\sigma(u)$
corresponding to the fractional Brownian motion, that is
\[
d\sigma(u)= \frac{H(1-2H)}{\Gamma(2-2H)\cos(\pi H)}
\left|u\right|^{1-2H}du, \quad H\in(0,1).
\]
Still for $n=1$, this set was further used in \cite{aal2} when
$d\sigma(u)$ is absolutely continuous with respect to Lebesgue
measure, with appropriate conditions on the Radon-Nikodym
derivative. Certain singular measures $d\sigma$ have been
considered in
\cite{MR2793121}.\\

We consider the hermitian form
\[
L_\sigma(\psi_1,\psi_2)=\int_{\mathbb
R^n}\widehat{\psi_1}(u)\overline{\widehat{\psi_2}(u)}d\sigma(u)
=\langle\widehat{\psi_1},\widehat{\psi_2}\rangle_{\mathbf
L_2(\sigma)},
\] and denote
\[
q_\sigma(\psi)=L_\sigma(\psi,\psi)
\]
the associated quadratic form. Such forms appear in particular in
the theory of generalized stochastic processes and generalized
stochastic fields. See \cite[Th\'eor\`eme $1^\prime$, p. 258, p.
283]{MR35:7123}. We now give another characterization of
$\mathcal C$ in Proposition \ref{pnfrechet} below. We first need
a preliminary lemma. In the statement of the lemma, $U_t$
denotes, for $t\in\mathbb R^n$ and $f\in \mathbf L_2(\mathbb
R^n,dx)$, the translation operator
\begin{equation}
\label{u_t}
(U_tf)(x)=f(x-t).
\end{equation}

\begin{La}
A positive semi-definite quadratic form $q$ continuous on
$\mathcal S_{\mathbb R^n}$ is of the form \eqref{qsigma}
\begin{equation*}
q(\psi)=q_\sigma(\psi)=\int_{\mathbb
R^n}|\widehat{\psi}(u)|^2d\sigma(u)
\end{equation*}
for some $\sigma\in\mathcal C$ if and only if
\begin{equation}
\label{q_u_t} q(U_t\psi)=q(\psi),\quad\forall \psi\in\mathcal
S_{\mathbb R^n}.
\end{equation}
\end{La}

{\bf Proof:} It is clear that $q_\sigma$ is a positive quadratic
form which satisfies \eqref{q_u_t} since
\[
\widehat{U_tf\psi}(u)=e^{iut}\widehat{\psi}(u).
\]
To prove the converse statement,
let $(\psi_p)_{p\in\mathbb N_0}$ be a sequence of elements in
$\mathcal D_{\mathbb R^n}$ which converge to some element
$\psi\in \mathcal D_{\mathbb R^n}$. Then there is a compact set
$K\subset\mathbb R^n$ such that all the $\psi_p$ and $\psi$ have
support inside $K$ and $(\psi_p)_{n\in\mathbb N_0}$ as well as
all the sequences of partial derivative converge uniformly to
$\psi$ and to the corresponding partial derivative respectively.
Since $K$ is bounded, it follows that, for every $\alpha$ and
$\beta$ in $\mathbb N_0^n$, the sequence $(x^\alpha
\psi_p^{(\beta)}(x))_{n\in\mathbb N_0}$ converge uniformly to
$x^\alpha \psi^{(\beta)}(x)$ on $K$, and hence on $\mathbb R^n$.
Thus, $(\psi_p)_{n\in\mathbb N_0}$ converges to $\psi$ in the
topology of $\mathcal S_{\mathbb R^n}$. Since $q$ is assumed
continuous in the topology of $\mathcal S_{\mathbb R^n}$ we have
\[
\lim_{p\rightarrow\infty}q(\psi_p)=q(\psi).
\]
Hence, $q$ is also continuous as a quadratic form on $\mathcal
D_{\mathbb R^n}$ Using \cite[Th\'eor\`eme 6, p. 160]{MR35:7123}
we see that $q$ restricted to $\mathcal D_{\mathbb R^n}$ is of
the form $q_\sigma$ for some $\sigma\in\mathcal C$. By continuity
it is equal to $q_\sigma$ on $\mathcal S_{\mathbb R^n}$.
\mbox{}\qed\mbox{}\\

\begin{Pn}
\label{pnfrechet}
Let $\sigma$ be a positive measure on $\mathbb R^n$. Then, $\sigma\in
\mathcal C$ if and only if the map
\[
\psi\mapsto q_\sigma(\psi)=\int_{\mathbb R^n}|\widehat{\psi}(u)|^2d\sigma(u)
\]
is continuous in the Fr\'echet topology of $\mathcal S_{\mathbb
R^n}$.
\end{Pn}

{\bf Proof:} One direction is done as in the proof of
\cite[Theorem 5.2]{MR2793121}: If $\sigma\in\mathcal C$ we can
write:
\[
\begin{split}
q_\sigma(\psi)&=\int_{\mathbb
R^n}|\widehat{\psi}(u)|^2(1+|u|^2)^p\frac{d\sigma(u)}{(1+|u|^2)^p}\\
&\le C\max_{u\in\mathbb
R^n}|\left(|\widehat{\psi}(u)|^2(1+|u|^2)^p\right),
\end{split}
\]
with $C=\int_{\mathbb R^n}\frac{d\sigma(u)}{(1+|u|^2)^p}$. We have
\[
\begin{split}
|\widehat{\psi}(u)|^2&=|\int_{\mathbb
R^n}\psi(x)\star\psi^\sharp(x)e^{-iu.x}dx|\\
&\le\int_{\mathbb R^n}
|\psi(x)\star\psi^\sharp(x)|dx\\
&\le\left(\int_{\mathbb R^n}|\psi|(x)dx\right)\left(\int_{\mathbb
R^n}|\psi^\sharp|(x)dx\right)\\
&=\left(\int_{\mathbb R^n}|\psi|(x)dx\right)^2.
\end{split}
\]
The terms involving powers of the components $u_j$ are treated in
the same  way.\\

Conversely, if $q_\sigma$ is continuous in the Fr\'echet
topology, then as in the preceding lemma, it is continuous as a
map from $\mathcal D_{\mathbb R^n}$ into $\mathbb C$. So
$q_\sigma$ is a translation invariant continuous sesquilinear
form from $\mathcal D_{\mathbb R^n}$ into $\mathbb C$. By the
previous lemma, there is a measure $\sigma_1\in\mathcal C$ such
that
\[
q_\sigma(\psi)=\int_{\mathbb
R^n}|\widehat{\psi}(u)|^2d\sigma_1(u),\quad\forall \psi\in
\mathcal D_{\mathbb R^n},
\]
and so $\sigma=\sigma_1$.
\mbox{}\qed
\mbox{}\\

In Theorem \ref{tm_Q} below we prove that there exists a bounded
linear operator $Q$ from $\mathcal S_{\mathbb R^n}$ into $\mathbf
L_2(d\sigma)$ such that
\begin{equation}
\label{L_Q} \int_{\mathbb
R^n}\widehat{\psi_1}(u)\overline{\widehat{\psi_2}}(u)d\sigma(u)
=\langle Q\psi_1,Q\psi_2\rangle_{\mathbf L_2(d\sigma)},\quad
\forall \psi_1,\psi_2\in\mathcal S_{\mathbb R^n}.
\end{equation}
The analysis of the problem uses Schwartz' kernel theorem and a
factorization result of Gorniak and Weron; see \cite{MR647140}
and \cite{MR750756}.
 By Schwartz' kernel
theorem there is a continuous linear positive operator from
$\mathcal S_{\mathbb R^n}$ into  $\mathcal S^\prime_{\mathbb
R^n}$ such that
\[
\int_{\mathbb
R^n}\widehat{\psi_1}(u)\overline{\widehat{\psi_2}(u)}d\sigma(u)=
\langle T\psi_1,\overline{\psi_2\rangle}.
\]
The form is positive, and so $T$ is a positive operator from
$\mathcal S_{\mathbb R^n}$ into  $\mathcal S^\prime_{\mathbb
R^n}$. It is proved in \cite{MR647140} that the space $\mathcal
S_{\mathbb R^n}$ has the factorization property, meaning that any
continuous positive operator from $\mathcal S_{\mathbb R^n}$ into
$\mathcal S^\prime_{\mathbb R^n}$ can be factorized via a Hilbert
space: There exists a Hilbert space $\mathcal H$ and a continuous
operator $Q$ from  $\mathcal S_{\mathbb R^n}$ into $\mathcal H$
such that $T=Q^*Q$. To conclude, it suffices to take an
isomorphism between $\mathcal H$ and $\mathbf L_2(\mathbb
R^n,dx)$. Note that the operator $Q$ will not, in general, be
bounded from $\mathbf L_2(\mathbb R^n,dx)$ into $\mathbf
L_2(d\sigma)$.\\

The argument below does not give an explicit construction. However, in the
proof of Theorem \ref{tm_Q} the operator $Q$ (denoted there by
$Q_\sigma$) is constructed in an explicit way.

\begin{Tm}
\label{tm_Q}
 Let $\sigma$ be a positive measure on $\mathbb R^n$
subject to \eqref{eq:bound}, and assume that
\[
{\rm dim}~\mathbf L_2(d\sigma)=\infty.
\]
There exists a continuous linear operator $Q_\sigma$ from $\mathcal
S_{\mathbb R^n}$ into $\mathbf L_2(d\sigma)$ such that
\eqref{L_Q} holds.
\end{Tm}

{\bf Proof:} We first note that the operators  $M_{u_k}$ of
multiplication by the variable $u_k$ in $\mathbf L_2(d\sigma)$ are
self-adjoint and commute with each other. Fix an isometric
isomorphism $W$ from $\mathbf L_2(d\sigma)$ onto ${\mathbf
L_2(\mathbb R^n,dx)}$, let $h\in\mathbf L_2(d\sigma)$ be defined
by
\[
h(u)=\frac{1}{(1+|u|^2)^{p/2}},
\]
and introduce $T_k=WM_{u_k}W^*$.
Define $T=(T_1,\ldots, T_n)$,
and
\begin{equation}
\label{eq:Q}
Q_\sigma\psi=(I+\sum_{k=1}^nT_k^2)^{p/2}\widehat{\psi}(T)Wh.
\end{equation}
Then, $Q_\sigma$ satisfies \eqref{L_Q}, as is seen using the
functional calculus for commuting normal operators. More
precisely, with
\[
W^*TW=(W^*T_1W,W^*T_2W,\ldots, W^*T_nW)
\]
in the third line, and
\[
M_u=(M_{u_1},M_{u_2},\ldots ,M_{u_n})
\]
in the fourth line, we have:
\[
\begin{split}
\|Q_\sigma\psi\|^2_{\mathbf L_2(\mathbb
R,dx)}&=\|(I+\sum_{k=1}^nT_k^2)^{p/2}\widehat{\psi}(T)Wh\|^2_{\mathbf
L_2(\mathbb R,dx)}\\
&=\|W^*(I+\sum_{k=1}^nT_k^2)^{p/2}\widehat{\psi}(T)Wh\|^2_{\mathbf
L_2(\sigma)}\\
&=\|(I+\sum_{k=1}^n(W^*T_kW)^2)^{p/2}\widehat{\psi}(W^*TW)h\|^2_{\mathbf
L_2(\sigma)}\\
&=\|(I+\sum_{k=1}^nM_{u_k}^2)^{p/2}\widehat{\psi}(M_u)h\|^2_{\mathbf
L_2(\sigma)}\\
&=\int_{\mathbb
R^n}\big|(1+|u|^2)^{p/2}\widehat{\psi}(u)\frac{1}{(1+|u|^2)^{p/2}}\big|^2d\sigma(u)\\
&=\int_{\mathbb R^n}|\widehat{\psi}(u)|^2d\sigma(u).
\end{split}
\]
\mbox{}\qed\mbox{}\\

\begin{La}
\label{la:dense}
The space $\mathcal S_{\mathbb R^n}$ is dense in
$\mathbf L_2(\sigma)$
\end{La}

{\bf Proof:} The space $\mathcal S_{\mathbb R^n}$ is closed under
pointwise product, and under conjugation. Furthermore, $\mathcal
S_{\mathbb R^n}$ separates points. Therefore, it is dense in
$C_0(\stackrel{\cdot}{\mathbb R^n})$, where
$\stackrel{\cdot}{\mathbb R^n}$ denotes the one point
compactification of $\mathbb R^n$, and where
$C_0(\stackrel{\cdot}{\mathbb R^n})$ denotes the space of
uniformly continuous functions $f$ defined on $\mathbb R^n$ and
such that $\lim_{|x|\rightarrow\infty}f(x)=0$.
\mbox{}\qed\mbox{}\\

\begin{Cy}
The completion of $\mathcal S_{\mathbb R^n}$ in the norm
$\left(\int_{\mathbb
R^n}|\widehat{\phi}(u)|^2d\sigma(u)\right)^{1/2}$ is $\mathbf
L_2(\mathbb R^n,dx)$.
\end{Cy}

{\bf Proof:} We have
\[
Q\phi=(I+T^2)^{p/2}\widehat{\phi}(T)W\left(\frac{1}{(1+|u|^2)^{p/2}}\right),
\]
and so
\[
\|Q\phi\|_{\mathbf L_2(\mathbb
R^n,dx)}=\|\widehat{\phi}\|_{\mathbf L_2(\sigma)}.
\]
By construction $W(\mathbf L_2(\sigma))=\mathbf L_2(\mathbb
R^n,dx)$, and this concludes the proof.
\mbox{}\qed\mbox{}\\

\begin{Tm}
Let $d\sigma\in\mathcal C$ be with the following property: For
every compact subset $K$ of $\mathbb R^n$, the exponential
functions $\left\{e_t(u)=e^{it\cdot u},\, t\in\mathbb
R^n\right\}$ are dense in $\mathbf L_2(K, d\sigma)$. Then, the
operator $Q_\sigma$ defined by \eqref{eq:Q} has dense range in
$\mathbf L_2(\mathbb R^n,dx)$.
\end{Tm}
{\bf Proof:} Let $g\in\mathbf L_2(\mathbb R^n,dx)$ be such that
\[
\langle Q_\sigma\psi,g\rangle_{\mathbf L_2(\mathbb
R^n,dx)}=0,\quad\forall\psi\in\mathcal S_{\mathbb R^n}.
\]
Thus
\[
\begin{split}
0&=\int_{\mathbb
R^n}\overline{(W^*g)(u)}\widehat{\psi}(u)d\sigma(u)\\
&=\iint_{\mathbb R^n\times\mathbb R^n} (W^*g)(u)e^{iu\cdot x}
\psi(x)d\sigma(u)dx\\
\end{split}
\]
Take $\psi$ which has support inside a fix compact $K$. It follows
that
\[
\int_{K}(W^*g)(u)e^{-iu\cdot x}d\sigma(u)=0\quad a.e.
\]
From the hypothesis we get that $(W^*g)(u)=0$ almost everywhere
on $K$. Since $K$ is arbitrary, we have that $W^*g=0$, and hence
$g=0$.\mbox{}\qed\mbox{}\\

\section{Closability}
\setcounter{equation}{0}
\label{sec:4}

As noted, our processes are indexed by a family $\mathcal C$ of
positive Radon measures $\sigma$ on $\mathbb R^n$, see Definition
3.1. The distinction between the two cases when $\sigma$ in
$\mathcal C$ is absolutely continuous with respect to Lebesgue
measure or not is crucial in applications of our processes to
stochastic integrals and Ito formulas.
 In this section we show that the absolute continuity
condition is equivalent to closability of a certain quadratic
form in $\mathbf L_2(\mathbb R^n)$.\\

 There are a variety of uses of
sesquilinear forms $q$ in operator theory and mathematical
physics, and we refer to the book \cite{Kato} for details. One
property for a sesquilinear form $q$ is "closability." This notion
depends on the choice of the ambient Hilbert space $\mathcal H$.
The notion of closability for sesquilinear forms plays a crucial
role in numerous applications. In the discussion below, the
ambient Hilbert space $\mathcal H$ will be $\mathbf L_2(\mathbb
R^n,
dx)$.\\

We recall the following definition.
\begin{Dn}
\label{dn:closable} The quadratic form $q$ defined on the
Schwartz space $\mathcal S_{\mathbb R^n}$ is closable if the
following condition holds: Given any sequence $(s_k)_{k\in\mathbb
N}$ of elements of $\mathcal S_{\mathbb R^n}$ for which
\[
\lim_{k\rightarrow\infty}\|s_k\|_{\mathbf L_2(\mathbb R^n,dx)}=0
\quad{and}\quad \lim_{k,\ell\rightarrow\infty}q(s_k-s_\ell)=0,
\]
it holds that
\[
\lim_{k\rightarrow\infty} q(s_k)=0.
\]
\end{Dn}

\begin{Tm}
\label{tm:closable}
The quadratic form $q_\sigma(\psi)$ is
closable if and only if $d\sigma$ is absolutely continuous with
respect to the Lebesgue measure.
\end{Tm}

{\bf Proof:} We first assume that $q$ is closable in the sense of
Definition \ref{dn:closable}, and we show that $d\sigma$ is
absolutely continuous with respect to the Lebesgue measure. We
divide the arguments into a number of steps.\\

STEP 1: {\sl There exists a positive selfadjoint linear operator
$\Gamma$ with domain containing $\mathcal S_{\mathbb R^n}$ and
such that
\begin{equation}
\label{q_sigma_kato}
q_\sigma(\psi)=\|\Gamma^{1/2}\psi\|^2_{\mathbf L_2(\mathbb
R^n,dx)},\quad \forall \psi\in S_{\mathbb R^n}.
\end{equation}}

This follows from Kato's theorem in \cite{Jo79,Kato,Ka66} for
closable quadratic forms. For the next step, recall that, for
$t\in\mathbb R^n$,
the translation operator $U_t$ has been defined in \eqref{u_t}.\\

STEP 2: {\sl The operator $\Gamma$ in \eqref{q_sigma_kato}
commutes with the translation operators $U_t$ in $\mathbf
L_2(\mathbb R^n,dx)$.}\\

Indeed, we have for $\psi\in\mathcal S_{\mathbb R^n}$
\[
\widehat{U_t\psi}(u)=e^{it\cdot u}\widehat{\psi}(u),\quad t\in\mathbb R^n,
\]
and so
\[
q_\sigma(\psi)=q_\sigma(U_t\psi),\quad\forall \psi\in\mathcal
S_{\mathbb R^n},
\]
that is
\[
\|\Gamma^{1/2}\psi\|^2_{\mathbf L_2(\mathbb R^n,dx)}=
\|\Gamma^{1/2}U_t\psi\|^2_{\mathbf L_2(\mathbb R^n,dx)}.
\]
By the uniqueness of the positive operator $\Gamma$ in
\eqref{q_sigma_kato} we obtain
\[
U_t^*\Gamma^{1/2} U_t=\Gamma^{1/2},
\]
and hence
\[
U_t^* \Gamma U_t=\Gamma.
\]
Hence the selfadjoint operator $\Gamma$ is in the commutant of
the unitary $n$-parameter group $\left\{U_t\right\}_{t\in\mathbb
R^n}$ acting on $\mathbf L_2(\mathbb R^n,dx)$.\\

STEP 3: {\sl $\Gamma$ is a convolution operator, that is there is
$m\in \mathbf L_1^{\rm loc}(\mathbb R^n,dx)$ such that
\begin{equation}
\Gamma\psi=m\star\psi.
\end{equation}
}\\

This is because the group $(U_t)_{t\in\mathbb R^n}$ is
multiplicity-free, and thus its commutant consists of convolution
operators. In our applications of closable quadratic forms we
make use of Kato's theory as presented in \cite{Ka66, Ka81} . The
thrust of Kato's theorem is that there is a precise way to
associate a selfadjoint operator to a closable quadratic form
defined on a dense subspace in a Hilbert space. We further make
use of results regarding unbounded operators commuting with
algebras of bounded operators. These results are in \cite{Jo79,
Jo84}, and \cite{Sto51}.\\

We can now conclude the proof:
\[
q_\sigma(\psi)=\int_{\mathbb R^n} (m\star
\psi)(x)\overline{\psi(x)}dx.
\]
Applying Parseval's equality on the right side of the above
equality we get
\[
q(\psi)=\int_{\mathbb R^n}\widehat{m}(u)|\widehat{\psi}(u)|^2du.
\]
We obtain thus that $\widehat{m}(u)\ge 0$ and
$d\sigma(u)=\widehat{m}(u)du$. It remains to prove that the
converse statement holds: if $d\sigma$ is absolutely continuous
with respect to Lebesgue measure, then $q_\sigma$ is closable.

\mbox{}\qed\mbox{}\\

\section{The generalized process $X_\sigma$}
\setcounter{equation}{0}

The focus of our paper is a family $\mathcal C$ of regular Borel
measures $\sigma$ on $\mathbb R^n$ (see Definition \ref{dn:C}).
We will be assigning a stationary-increment processes to every
$\sigma$ in the set $\mathcal C$ (see Theorem \ref{thm:central}). The
properties of these processes will be studied in the rest of the
paper. In Section \ref{sec:8} we show that $\mathcal C$ is closed
under convolution of measures.

\begin{Dn}
We denote by $\mathbf L(\mathcal S_{\mathbb R^n},\mathbf
L_2(\mathbb R^n,dx))$ the space of linear operators $Q$ from
$\mathcal S_{\mathbb R^n}$ to $\mathbf L_2(\mathbb R^n,dx)$.
\end{Dn}

In the following statement, recall that we have set
$\Omega=\mathcal S_{\mathbb R^n}^\prime(\mathbb R)$.

\begin{Tm}
\label{thm:central}
Let $\sigma\in\mathcal C$. Then, there exists
a probability measure $\mu_\sigma$ on $\mathcal S^\prime_{\mathbb
R^n}(\mathbb R)$, an element $Q_\sigma\in \mathbf L(\mathcal
S_{\mathbb R^n},\mathbf L_2(\mathbb R^n,dx))$ and a generalized
Gaussian stochastic process
$\left\{X_\sigma(\psi)\right\}_{\psi\in\mathcal S_{\mathbb R^n}}$
such that:
\begin{align}
\label{1}
X_\sigma(\psi)&=\widetilde{(Q_\sigma(\psi))},\quad\forall
\psi\in\mathcal S_{\mathbb R^n}\\
\intertext{and} \int_{\Omega}e^{i\langle \w,
X_\sigma(\psi)\rangle}d\mu_W(\w)&= \int_{\Omega}e^{i\langle
\w,\psi\rangle}d\mu_\sigma(\w) \label{2}
=e^{-\frac{\|\widehat{\psi}\|^2_{\mathbf L_2(\sigma)}}{2}}.
\end{align}
\end{Tm}

{\bf Proof:} Using Theorem \ref{tm_Q}, we take an operator
$Q_\sigma\in \mathbf L(\mathcal S_{\mathbb R^n},\mathbf
L_2(\mathbb R^n,dx))$ such that \eqref{L_Q} holds. We define a
generalized stochastic process $(X_\sigma)_{\psi\in\mathcal
S_{\mathbb R^n}}$ via the formula \eqref{1}. By definition of
$d\mu_W$ we have
\begin{equation}
\label{eq:Qsigma} \int_{\Omega} e^{i\langle \w,
X_\sigma(\psi)\rangle}d\mu_W(\w)=e^{-\frac{\|Q_\sigma\psi\|^2_{\mathbf
L_2(\mathbb R^n,dx)}}{2}}.
\end{equation}

\mbox{}\qed\\\

The operator $Q_\sigma$ is continuous from $\mathcal S_{\mathbb
R^n}$ into $\mathbf L_2(\mathbb R^n,dx)$. Its adjoint is a
continuous operator from $\mathbf L_2(\mathbb R^n,dx)$ into
$\mathcal S^\prime_{\mathbb R^n}$. We have the diagram

\[\renewcommand{\arraystretch}{1.5}
\begin{array}{ccccc}
\mathcal S_{\mathbb R^n}&\stackrel{i}{\longhookrightarrow}
&\mathbf L_2(\mathbb R^n,dx)& \stackrel{i^*}{\longhookrightarrow}
&\mathcal S_{\mathbb
R^n}^\prime\\
&\stackrel{Q_\sigma}{\searrow}& &\stackrel{Q_\sigma^*}{\nearrow}\\
\mathcal S_{\mathbb R^n}&\stackrel{i}{\longhookrightarrow}
&\mathbf L_2(\mathbb
R^n,dx)&\stackrel{i^*}{\longhookrightarrow}&\mathcal S_{\mathbb
R^n}^\prime
\end{array}.\]

\begin{Tm}
The operator
\[
Q_\sigma^*Q_\sigma\,\,\,:\,\,\, \mathcal S_{\mathbb
R^n}\longrightarrow \mathcal S^\prime_{\mathbb R^n}
\]
is given by the formula
\[
(Q_\sigma^*Q_\sigma\psi)(x)=\widehat{\sigma}(x-\cdot)(\psi),
\]
meaning that
\begin{equation}
\label{eq:Q*Q}
\langle Q_\sigma^*Q_\sigma \psi,\phi\rangle=\int_{\mathbb
R^n}d\sigma(u)\widehat{\psi}(u)\overline{\widehat{\phi}(u)}du.
\end{equation}
\end{Tm}

{\bf Proof:} By definition of $\widehat{\sigma}$ we have
\[
\langle \widehat{\sigma}(t-s),v(t,s)\rangle=\int_{\mathbb R^n}\left(
\int_{\mathbb
R^n\times\mathbb R^n} e^{-iu(t-s)}v(t,s)dtds\right)d\sigma(u).
\]
Taking $v(t,s)=\psi(t)\overline{\phi}(s)$ we obtain
\[
\langle \widehat{\sigma}(t-s),\psi(t)\phi(s)\rangle=\int_{\mathbb
R^n}\psi(u)\overline{\phi(u)}\sigma(u).
\]
But, by definition of the adjoint we have
\[
\langle Q_\sigma\psi,Q_\sigma\phi\rangle=\langle
Q_\sigma^*Q_\sigma \psi, \phi\rangle
\]
and hence the result.
\mbox{}\qed\mbox{}\\

For related computation when $n=p=1$ see \cite[Section 2]{aal2}.
\section{Reproducing kernels}
\setcounter{equation}{0}
In this section we study the reproducing kernel Hilbert space
$\mathscr H(C_\sigma)$ associated to the positive definite
function
\[
C_\sigma(\phi-\psi)=e^{-\frac{\|\widehat{\phi}-\widehat{\psi}\|^2_{\mathbf
L_2(\sigma)}}{2}},
\]
where $\phi$ and $\psi$ vary through the space
$\mathcal S_{\mathbb R^n}(\mathbb R)$
of {\sl real-valued} Schwartz functions.
Let us already mention that the
fact that the functions are real play a key role in the arguments.
We set for $f\in\mathbf L_2(\mu_\sigma)$
\begin{equation}
\label{eq:f}
(\mathscr F_\sigma f)(\psi)=\int_\Omega e^{-i\langle
\w,\psi\rangle}f(\w)d\mu_\sigma(\w).
\end{equation}
\begin{Tm}
\label{tm:iso}
The map $\mathscr F_\sigma$ is isometric from
$\mathbf L_2(\mu_\sigma)$ onto $\mathscr H(C_\sigma)$.
\end{Tm}
{\bf Proof:} We set
\[
e_\phi(\w)=e^{i\langle \w,\phi\rangle}.
\]
Since we are in the real case $e_\phi$ has modulus $1$ and hence
belongs to $\mathbf L_2(\mu_\sigma)$. We note that
\[
(\mathscr F_\sigma e_\phi)(\psi)= \int_\Omega e^{-i\langle
\w,\psi\rangle}e^{i\langle \w,\phi\rangle}d\mu_\sigma(\w)=
C_\sigma(\phi-\psi).
\]
It follows from this equality that, for $f$ in the closed
linear span of the functions $e_\psi$ in
$\mathbf L_2(\Omega, d\mu_\sigma)$ the function \eqref{eq:f}
belongs to $\mathscr H(C_\sigma)$ and that
\[
\|\mathscr F_\sigma f\|_{\mathscr H(C_\sigma)}=\|f\|_{\mathbf L_2(\Omega,
d\mu_\sigma)}.
\]
To conclude the proof of the theorem we need to show that the
above closed linear span is equal to $\mathbf L_2(\Omega,
d\mu_\sigma)$. The argument is a bit long, and divided into a
number of steps.\\

STEP 1: {\sl The closed linear span of the $e_\psi$ is dense in
$\mathbf L_2(\Omega, d\mu_\sigma)$ if and only if the polynomials
are dense in $\mathbf L_2(\Omega, d\mu_\sigma)$.}\\

Indeed, one direction is clear from the power series expansion of
the exponential. To prove the converse, consider $\phi$ of the
form
\[
\phi=\sum_{m=1}^M \lambda_m\phi_m,
\]
where $M\in\mathbb N$, $\lambda_1,\ldots, \lambda_M\in\mathbb R$
and $\phi_1,\ldots, \phi_M\in\mathcal S_{\mathbb R^n}(\mathbb R)$, and let
$F\in\mathbf L_2(\Omega,\mu_\sigma)$ be
orthogonal to all the $e_\phi$.
\[
\int_{\Omega}F(\w)e_\phi(\w)d\mu_\sigma(\w)= \int_{\Omega}F(\w)e^{\lambda_1\langle
\w,\phi_1\rangle}\cdots e^{\lambda_1\langle
\w,\phi_M\rangle}d\mu_\sigma(\w).
\]
Differentiating this expression with respect to $\lambda_1$ and
setting $\lambda_m=0$, $m=1,\ldots M$ we get
\[
\int_{\Omega}F(\w)\langle
\w,\phi_1\rangle d\mu_\sigma(\w).
\]
More generally, differentiating with respect to all the
variables, we get orthogonality with respect to all the
polynomials.\\

We denote by $\mathcal H$ the closed linear span of the random
variable $\w\mapsto\langle \w,\phi\rangle$ in $\mathbf
L_2(\Omega, \mu_\sigma)$.\\

STEP 2: {\sl The map $\mathscr F_\sigma$ which to the random
variable $\w\mapsto\langle \w,\phi\rangle$ associates the function
$\widehat{\phi}$ extends to a unitary from $\mathcal H$ onto
$\mathbf L_2(\sigma)$.}\\

Isometry is is a direct consequence of \eqref{2}. The fact that
$\mathscr F_\sigma$ is onto follows from Lemma \ref{la:dense}.\\

We note that $\Gamma(\mathcal H)=\mathbf L_2(\Omega,
\mu_\sigma)$, where $\Gamma$ denotes the second quantization
functor. We denote by $\Gamma(\mathscr F_\sigma)$ the extension of
this map from $\mathbf L_2(\Omega, \mu_\sigma)$ into
the Fock space $\mathcal F(\mathbf L_2(\sigma))$.\\

In the next step, we denote by $H_m$ and $h_m$ the Hermite
polynomials and the Hermite functions respectively,
$m=0,1,\ldots$. Furthermore, we set $\ell$ to be the set of
sequences
\begin{equation}
\label{eq:ell}
(\alpha_1,\alpha_2,\ldots),
\end{equation}
indexed by ${\mathbb N}$  with values in ${\mathbb N}_0$, for
which only a finite number of elements $\alpha_j\not=0$.\\

STEP 3: {\sl The functions $ H_\alpha(\w)=\prod_{m=1}^\infty H_{\alpha_{m}}(\langle
\w,h_{m}\rangle)$ form an orthogonal system in $\mathbf
L_2(\Omega, \mu_\sigma)$ when $\alpha$ runs through $\ell$.}\\

See for instance \cite[Theorem 22.24, p. 26]{new_sde}.\\

STEP 4: {\sl The span is dense}.\\

It is enough to prove that the polynomials are dense in $\mathbf
L_2(\Omega, \mu_\sigma)$. This is turn
follows from the fact that the "exponentials"
\[
\Gamma(\phi)=\sum_{m=1}^\infty
\frac{1}{\sqrt{m!}}\phi\otimes\cdots\otimes\phi
\]
are dense in $\Gamma(\mathbf L_2(\sigma))$. See for instance
\cite[Proposition 2.14, p. 28]{MR1851117} for the latter.

\mbox{}\qed\mbox{}\\

\section{The associated stochastic process, second construction,
and the fundamental isomorphism}
\setcounter{equation}{0}
\label{sec:7}
To understand the stochastic processes governed by measures
$\sigma$ from $\mathcal C$, we will be making use of a
fundamental isomorphism for the associated white noise spaces. We
turn to the details below. We denote by $\mathcal W_\sigma$ the
space
\[
\mathcal W_\sigma=\mathbf L_2(\mathcal S^\prime_{\mathbb
R^n}(\mathbb R),\mathcal B(\mathscr S^\prime_{\mathbb
R^n}(\mathbb R)),d\mu_\sigma).
\]
We present another Gaussian process with generalized covariance
$q_\sigma(\psi_1,\psi_2)$ and an isomorphism between $\mathcal W$
and a closed subspace of $\mathcal W_\sigma$. Equation \eqref{2}
can be rewritten as
\begin{equation}
\label{eqY}
E[e^{iY_\sigma(\psi)}]=e^{-\frac{\|\widehat{\psi}\|^2_{\mathbf
L_2(d\sigma)}}{2}},
\end{equation}
where $Y_\sigma(\psi)$ denotes the random variable defined by
\[
Y_\sigma(\psi)=\langle \omega, \psi \rangle_{\mathcal
S^\prime_{\mathbb R^n}(\mathbb R),\mathcal S_{\mathbb
R^n}(\mathbb R)}.
\]
We have now the generalization of the map \eqref{eq:iso} for
$\sigma\in\mathcal C$. We will denote this map
$\widetilde{\psi}^\sigma$, so that
\begin{equation}
\label{tildesigma}
\widetilde{\psi}^\sigma=Y_\sigma(\psi).
\end{equation}
A first consequence of the fact that \eqref{tildesigma} is an
isometry is:

\begin{Tm}
\label{tm:iso1}
Let $\psi, \psi_1, \psi_2\in\mathcal S_{\mathbb R^n}(\mathbb R)$.
Then:
\begin{eqnarray}
E[Y_\sigma(\psi)]&=&0, \\
E[Y_\sigma(\psi_1)\overline{Y_\sigma(\psi_1)}]&=&\int_{\mathbb
R^n} \widehat{\psi_1}(u)\overline{\widehat{\psi_2}(u)}d\sigma(u).
\end{eqnarray}
\end{Tm}

\begin{Tm}
The map $\psi\mapsto Y_\sigma(\psi)$ extends to an isometric map
from $\mathbf L_2(\mathbb R^n,dx)$ into $\mathcal W_\sigma$.
\end{Tm}

{\bf Proof:} In \eqref{eqY} replace $\psi$ by $\epsilon\psi$ with
$\epsilon\in\mathbb R$. Developing along the powers of $\epsilon$
both sides of \eqref{eqY} leads then to the result.
\mbox{}\qed\mbox{}\\

\section{Ergodicity}
\setcounter{equation}{0}
\label{sec:ergodicity}
%
%

%
%
We start with Theorem \ref{thm:central}. Apply \eqref{2} to the
function $\psi_t$ defined by
\[
\psi_t(x)=\psi(t+x),
\]
we obtain
\[
\int_{\mathbb R^n}e^{i\langle \w,
X_\sigma(\psi_t)\rangle}=e^{-\frac{1}{2}\|\widehat{\psi_t}
\|^2_{\mathbf L_2(\sigma)}}=e^{-\frac{1}{2}\|\widehat
{\psi}\|^2_{\mathbf L_2(\sigma)}}
\]
since
\[
\widehat{\psi_t}(x)=e^{itx}\widehat{\psi}(x).
\]
Therefore, we have that, for every $\psi\in\mathscr S_{\mathbb
R^n}(\mathbb R)$ and every $t\in\mathbb R$,
\[
\int_{\mathbb R^n}e^{i\langle \w,
X_\sigma(\psi_t)\rangle}d\mu_W(\w)= \int_{\mathbb R^n}e^{i\langle
\w, X_\sigma(\psi)\rangle}d\mu_W(\w).
\]
For every $t\in\mathbb R^n$ the generalized processes
\[\{X_\sigma(\psi)\}_{\psi\in\mathscr S_{\mathbb R^n}(\mathbb
R)}\quad {\rm and}\quad \{X_\sigma(\psi_t)\}_{\psi\in\mathscr
S_{\mathbb R^n}(\mathbb R)}
\]
have the same generating function. It follows that there is a
unitary map $U_t$ such that
\[
U_tX(\psi)=X(\psi_t),\quad \psi\in\mathcal S_{\mathbb R^n}(\mathbb
R).
\]
The family $(U_t)_{t\in\mathbb R^n}$ clearly forms a group.

\begin{La}
Let $\mathcal H$ be a Hilbert space and $t\mapsto X(t)$ a
function from $\mathbb R^n$ into $\mathcal H$. Then, the
following are equivalent:\\
$(1)$ There exists a function $F_X$ on $\mathbb R^n$ such that
\begin{equation}
\label{eq:eq1}
\|X(t)-X(s)\|_{\mathcal H}^2=F_X(t-s),\quad\forall t,s\in\mathbb
R^n.
\end{equation}
$(2)$ For every $t\in\mathbb R^n$, the map $U(t)$ defined by
\begin{equation}
U(t)X(s)=X(t+s),\quad s\in\mathbb R^n,
\label{eq:eq2}
\end{equation}
is unitary and satisfies
\begin{equation}
\label{eq:eq3}
U(t_1+t_2)=U(t_1)U(t_2),\quad \forall t_1,t_2\in\mathbb R^n.
\end{equation}
\label{la1}
\end{La}
{\bf Proof:} Assume that \eqref{eq:eq1} holds. Then define $U(t)$
on the span of the linear span of the elements
$\left\{X(s)\,\,;,\,\ s\in\mathbb R^n\right\}$ by \eqref{eq:eq2}.
Using \eqref{eq:eq1}, one checks that $U(t)$ extends from the
subspace ${\rm span}~\left\{X(s)\,\,;,\,\ s\in\mathbb R^n\right\}$
to all of $\mathcal H$, and that the extension, still denoted by
$U(t)$, satisfies \eqref{eq:eq3}.\\

Conversely, assume that $(2)$ is in force. Then
\[
\begin{split}
\|X(t)-X(s)\|_{\mathcal H}^2&=\|U(t)X(0)-U(s)X(0)\|_{\mathcal
H}^2\\
&=\|U(t-s)X(0)-X(0)\|^2_{\mathcal H}, \quad
\end{split}
\]
using \eqref{eq:eq3} and the unitarity assumption, and thus
\eqref{eq:eq1} holds.
\mbox{}\qed\mbox{}\\

We denote by $U(\mathcal H)$ the set of unitary operators from
$\mathcal H$ into itself. The proof of the following lemma is
easy, and will be omitted.

\begin{La}
Let $(X(t))_{t\in\mathbb R^n}$ be a process as defined in the
preceding lemma, and assume \eqref{eq:eq1}. Let
\[
U\,\,:\,\ \mathbb R^n\longrightarrow U(\mathcal H)
\]
denote the corresponding representation of the group $(\mathbb
R^n,+)$ acting by unitary operators on $\mathcal H$. The, the
following are equivalent:\\
$(1)$ $(U(t))_{t\in\mathbb R^n}$ is a strongly continuous
representation.\\
$(2)$ The function $F_X$ defined by \eqref{eq:eq1} is continuous
from $\mathbb R^n$ into $\mathcal H$.
\label{la2}
\end{La}

We apply the preceding results to the special case $\mathcal
H=\mathbf L_2(\Omega,\mathcal B_\Omega,\mu_\sigma)$ (recall that
$\Omega=\mathcal S^\prime_{\mathbb R^n}(\mathbb R)$; see
\eqref{eq:Omega}) is the Wiener space corresponding to the
measure $\mu_\sigma$ as in Section \ref{sec:7}, and where
$Y_\sigma$ is the corresponding stationary-mean-square process.
The conditions in Lemmas \ref{la1} and \ref{la2} are satisfied,
and we say that $(U_\sigma(t))_{t\in\mathbb R^n}$ is the
corresponding unitary representation of the group $(\mathbb
R^n,+)$. Then:
\begin{equation}
\label{eq:eq6}
(U_\sigma(t)f)(\w)=f(\w(\cdot-t)),\quad\forall
\w\in\mathcal S_{\mathbb R^n}^\prime\,\,{\rm and}\,\,\forall f\in
\mathbf L_2(\Omega,\mathcal B_\Omega,\mu_\sigma).
\end{equation}

Indeed, by Lemma \ref{la1}, it is enough to determine
$U_\sigma(t)$ on the linear span of the elements $Y_\sigma(s)$
when $s$ runs through $\mathbb R^n$. For a given $s\in\mathbb
R^n$ and  $f=Y_\sigma(s)$ we have:
\begin{equation}
\label{eq:eq7}
U_\sigma(t)Y_\sigma(s)=Y(s+t),
\end{equation}
so that
\begin{equation}
\left(U_\sigma(t)Y_\sigma(s)\right)(\w)=(Y_\sigma(s+t))(\w).
\label{eq:eq8}
\end{equation}

We now compute \eqref{eq:eq6} for $f=Y\sigma(s)$. We have
\[
\begin{split}
\left(U_\sigma(t)Y_\sigma(s)\right)(\w)&=\left(Y_\sigma(s+t)\right)(\w)
\quad\mbox{(\rm where we have used \eqref{eq:eq7})}\\
&=\w(s+t)\\
&=Y_\sigma(\w(\cdot -t)),
\end{split}
\]
by \eqref{tildesigma}, which is the desired conclusion.\\

We conclude this section with:

\begin{Tm}
Given $\sigma\in\mathcal C$, let $X_\sigma$ be the corresponding generalized
stochastic process defined by \eqref{1},, and let $\mu_\sigma$ denote the measure
on $\Omega=\mathcal S^\prime_{\mathbb R^n}(\mathbb R)$
from Theorem \ref{thm:central}. Then the following two conditions
are equivalent:\\
$(i)$ The representation $(U_\sigma(t))_{t\in\mathbb R^n}$ in
\eqref{eq:eq7} has
$0$ as its only point spectrum, with a one-dimensional eigenspace.\\
$(ii)$ The action of $\mathbb R^n$ on $\Omega$:
\[
(t,\w)\mapsto\w(\cdot+t)
\]

is ergodic with respect to the measure $\mu_\sigma$
\end{Tm}

{\bf Proof:} It suffices to combine Theorem \ref{thm:central} and Lemmas
\ref{la1} and \ref{la2}.
\mbox{}\qed\mbox{}\\

\section{The Hilbert space of sigma-functions and applications}
\setcounter{equation}{0}
\label{sigma_function}

In Sections 3 through 7, we
     introduce a family of Gaussian stochastic processes as follows.
     Our construction goes beyond what has been done before in a number
     of ways, as we proceed to outline. The processes $Y$ we consider
     are indexed by the Schwartz space $\mathcal S_{\mathbb R^n}$
     (i.e., random variable-valued tempered distributions),
     and it is assumed that the expectation of the square of all the increments by
     $Y$ is independent, not the increments themselves; i.e., that for fixed $Y$,
     the expectation of the squares of differences  $Y(\psi)$  depends only
     on the differences of points $\psi$  in the space of test functions
     $\mathcal S_{\mathbb R^n}$. This allows us then to determine a given process $Y$ , via
     its covariance function, from a corresponding measure $d\sigma$  on
     $\mathbb R^n$. Conversely, we show that for every measure $d\sigma$ on
     $\mathbb R^n$, and falling  in a suitably defined class $\mathcal C$, there is
     a uniquely defined process  $Y_\sigma$  determined by its two-point covariance
     functions. Given $Y_\sigma$, the formula (see Theorem \ref{tm:iso1}) for its
     covariance is given by an integration with respect to $d\sigma$.
In this section, we will be interested in the correlations
computed from two such processes, each determined from a
different measure $d\sigma$   in the class $\mathcal C$.
  But this then introduces a difficulty as there are two measures, and
  they might be relatively singular, or perhaps they may allow for a
  suitable comparison, for example with the use of a Radon-Nikodym derivative.
So it is not at all clear how to compare and to compute the
covariance for the different processes. In order to get around
this we introduce a Hilbert space $\mathscr H$ of equivalence
classes, each equivalence class determined by some measure
$d\sigma$  from our class $\mathcal C$. More precisely, $\mathscr
H$ consists of equivalence classes of pairs $( \psi, d\sigma )$
where $\psi$ and $d\sigma$ are related, i.e., the function $\psi$
is assumed in $\mathbf L_2(d\sigma)$. We are then able to write
down an isometric transform for our processes which solves the
comparison problem, i.e., comparing two processes computed from
different measures $d\sigma$. Our transform is between the
Hilbert space $\mathscr H$ and the $\mathbf L_2$ space of Wiener
white noise measure on the space of real tempered Schwartz
distributions. We note that the Hilbert space $\mathscr H$ ,
often called a Hilbert space of sigma-functions, has been used in
the literature for a variety of different unrelated purpose, for
example spectral theory \cite{Nelson_flows}, infinite products of
measures \cite{Kak48}, harmonic
analysis \cite{Sch95} , probability theory \cite{Le62}, and more.\\

Let $\sigma_1$ and $\sigma_2$ be in $\mathcal C$, and let
$f_i\in\mathbf L_2(d\sigma_i)$ for $i=1,2$. Following
\cite{Nelson_flows}, we say that the two pairs $(f_1,\sigma_1)$
and $(f_2,\sigma_2)$ are equivalent if
\begin{equation}
\label{eq:radon}
f_1\sqrt{\frac{\rm d\sigma_1}{\rm d\lambda}}=
f_2\sqrt{\frac{\rm d\sigma_2}{\rm d\lambda}},\quad \lambda\,\
a.e.,
\end{equation}
where $\lambda$ is a measure such that both $\sigma_1$ and
 $\sigma_2$ are absolutely continuous with respect to $\lambda$,
 and where,
for instance, $\frac{\rm d\sigma_1}{\rm d\lambda}$ denotes the
Radon-Nikodym derivative. It is enough to check \eqref{eq:radon}
for $\lambda=\sigma_1+\sigma_2$. The equality \eqref{eq:radon}
defines indeed an equivalence relation. We denote by
$f\sqrt{d\sigma}$ the equivalence class of $(f,\sigma)$, and by
$\mathscr H$ the space of all such equivalent classes. The form
\[
\langle f_1\sqrt{d\sigma_1},f_2\sqrt{d\sigma_2}\rangle_{\mathscr
H}=\int_{\mathbb R^n}f_1(x)\overline{f_2(x)}\sqrt{\frac{\rm
d\sigma_1}{\rm d\lambda}\frac{\rm d\sigma_2}{\rm
d\lambda}}d\lambda
\]
is a well defined inner product, and the space $\mathscr H$
endowed with this inner product is a Hilbert space. This space
has been studied by numerous authors, and in particular by
Kakutani, L\'evy, Schwartz and Nelson. See \cite{Kak48},
\cite{MR0099595}, \cite{Nelson58}.\\

We now define the generalization of the operator $T_m$ in the
present setting by
\begin{equation}
\label{eq:Rsigma}
\widehat{R_\sigma f}=\widehat{f}\sqrt{\frac{\rm d\sigma}{\rm
d\lambda}}.
\end{equation}
With the use of the space $\mathscr H$ we now prove:

\begin{Tm}
Let $\sigma_1$ and $\sigma_2$ be two measures in $\mathcal C$,
and let $X_{\sigma_1}$ and $X_{\sigma_2}$ be the associated
generalized processes given by \eqref{eq:Qsigma}:
\[
\int_{\mathcal W} e^{i\langle
X_{\sigma_k}(\psi),\w\rangle}d\mu_W(\w)=
e^{-\frac{1}{2}\int_{\mathbb
R^n}|\widehat{\psi}(u)|^2d\sigma_k(u)},\quad k=1,2.
\]
Then, $X_{\sigma_1}$ and $X_{\sigma_2}$ are orthogonal in the
white  noise space $\mathcal W$ if and only if the measures
$\sigma_1$ and $\sigma_2$ are mutually singular.
\label{tm:9.1}
\end{Tm}

{\bf Proof:} We denote $W_\sigma$ the map which to $f\in\mathbf
L_2(\rm d\sigma)$ associates the equivalence class $f\sqrt{\rm
d\sigma}$, and by $\mathscr H(\sigma)=W_\sigma(\mathbf L_2(\rm
d\sigma))\subset\mathscr H$. Then the spaces $\mathscr
H(\sigma_1)$ and $\mathscr H(\sigma_2)$ are orthogonal in
$\mathscr H$ if and only if $\sigma_1$ and $\sigma_2$ are
mutually singular. Define operators $R_1$ and $R_2$ as in
\eqref{eq:Rsigma}, that is:
\[
\widehat{R_kf}=\widehat{f}\sqrt{\frac{\rm d\sigma_k}{\rm d\sigma
}},\quad k=1,2,
\]
with $\sigma=\sigma_1+\sigma_2$. We have
\[
\begin{split}
\langle \widehat{f}\sqrt{\rm d\sigma_1}, \widehat{f}\sqrt{\rm
d\sigma_2}\rangle_{\mathscr H}&=\int_{\mathbb R^n} \widehat{f}
\sqrt{\frac{\rm d\sigma_1}{\rm d\sigma }} \sqrt{\frac{\rm
d\sigma_2}{\rm d\sigma }}\overline{\widehat{g}}d\sigma\\
&=\langle X_\sigma(R_1f), X_\sigma(R_2g)\rangle_{\mathcal W}\\
&=\langle X_{\sigma_1}(f),X_{\sigma_2}(g)\rangle_{\mathcal W}.
\end{split}
\]
\mbox{}\qed\mbox{}\\

\section{Generalized Karhunen-Lo\`eve expansion}
\setcounter{equation}{0}
Following \cite{MR0345194} we define a random measure to be a
countably additive function on a given sigma algebra $\mathcal B$,
taking values in a space of random variables on a probability
space. If moreover disjoint sets in $\mathcal B$ are mapped into
independent random variables we say that the random measure is a
Wiener process. Starting with one of our stochastic processes $X$
as introduced in the first section, we show that it admits a
direct integral decomposition along an essentially unique Wiener
process $Z_X$, depending on $X$, and  with $Z_X$  Gaussian, i.e.,
Gaussian distributions in the fibers.\\

The classical Karhunen-Lo\`eve theorem (see e.g.,
\cite{ash_book}) applies to a restricted family of Gaussian
processes, i.e., to a system of stochastic processes
$\left\{X(t)\right\}$ indexed by  $t\in \mathbb R$ (typically to
represent time), and with specified joint distributions. When it
applies, it offers a separation of variables, expanding the
process $\left\{X(t)\right\}$ as a countable direct sum of a
system of independent identically distributed ( i.i.d.) $N(0,
1)$ random variables. Below (Theorem \ref{tm:7.2}) we extend this
to apply to the most general family of Gaussian processes (from
our section 5) . A key step in our more general expansion is a
systematic theory of direct integrals, taking the role of direct
sum-expansions in the restricted setting.\\

The theme of this section is the study of normal fields and their
role in direct integral decompositions of the stochastic
processes we introduced in sect 5 above. Our normal fields
(Definition \ref{def:normal}) extend a notion of Wiener processes
as introduced in \cite{MR545382}. In fact the study of Wiener
processes was initiated in a special case in early papers by
Ito.\\

 A random measure over a fixed sigma-algebra
$\mathcal M$ in a given measure space is a countably additive
mapping from $\mathcal M$ into random variables of some
probability space $(\Omega, P)$ , typically with $P$ some fixed
path-space measure. Alternatively, random measures are also known
as stochastic processes indexed by a measure space. If a random
measure takes stochastically independent values on disjoint sets
in $\mathcal M$, it is called a Wiener process. Indeed, Wiener
processes were extensively studied in \cite{MR666306},
\cite{MR545382}, see also \cite{Hida_BM}. Examples of Wiener
processes are Poisson processes, normal distributions, and jump
processes. A theorem in \cite{MR545382} states that every Wiener
process naturally decomposes into a sum of three components, a
Poisson, a normal, and jump process.

\begin{Dn}
Let $\mathcal B(\mathbb R^n)$ denote the sigma-algebra of Borel
sets of $\mathbb R^n$. Fix a measure $\sigma$ from $\mathcal C$. A
normal field will be a function
\[
Z\,\,:\,\,\,\mathcal B(\mathbb R^n)\times \Omega\longrightarrow
\mathbb R,
\]
with the following properties:
\begin{enumerate}
\item For every $\omega\in\Omega$, the function $Z(\cdot,\w)$ is
Borel-measurable on $\mathbb R^n$.

\item For every Borel set $A$, $Z(A,\cdot)\in\mathcal W_\sigma$ and is Gaussian.

\item $E_\sigma[Z(A,\cdot)]=0$.

\item For every $A_1,A_2\in\mathcal B(\mathbb R)$,
\[
E_\sigma[Z(A_1,\cdot)Z(A_2,\cdot)]=\int_{A_1\cap
A_2}\frac{d\sigma(u)}{(1+u^2)^p}.
\]
\end{enumerate}
\label{def:normal}
\end{Dn}

The following theorem can be seen as a generalized
Karhunen-Lo\`eve expansion for a special subfamily of elements of
$\mathcal C$ (namely, $n=p=1$).

\begin{Tm}
Let $d\sigma$ be a positive measure on $\mathbb R$ such
that
\[
\int_{\mathbb R}\frac{d\sigma(u)}{u^2+1}<\infty.
\]
Then there exists a normal field $Z(du,\cdot)$ such that
\begin{equation}
(X_\sigma(t))(\w)=\int_{\mathbb
R^n}\sqrt{1+|u|^2}\frac{e^{iu.t}-1}{u}Z(du,w),\quad \w\in\Omega.
\label{eq:normal}
\end{equation}
\label{tm:7.2}
\end{Tm}

{\bf Proof:} We proceed in a number of steps.\\

STEP 1: { \sl Construction of a projection-valued measure.}\\

As in the proof of Theorem \ref{tm_Q} we set $M_{u}$ denote the
operator of multiplication by the variable $u$ in $\mathbf
L_2(d\sigma)$ (recall that here $n=1$) and we fix an isometric
isomorphism $W$ from $\mathbf L_2(d\sigma)$ onto ${\mathbf
L_2(\mathbb R,dx)}$, and let
\begin{equation}
\label{eq:h}
h(u)=\frac{1}{\sqrt{1+u^2}}.
\end{equation}
The spectral theorem applied to the selfadjoint operator
$T=WM_uW^*$ leads to a projection-valued measure
\[
P\,\,\,\,:\,\,\,{\mathcal B}(\mathbb R)\longrightarrow {\rm
Proj}~(\mathbf L_2(\mathbb R,dx)),
\]
such that
\[
\begin{split}
T&=\int_{\mathbb R}\lambda P(d\lambda),\\
\int_{\mathbb R}\|P(d\lambda)f\|^2_{\mathbf L_2(\mathbb R,dx)}&=
\|f\|^2_{\mathbf L_2(\mathbb R,dx)},\quad \forall f\in\mathbf
L_2(\mathbb R,dx).
\end{split}
\]
Recall that one defines $x(T)=\int_{\mathbb R}
x(\lambda)P(d\lambda)$ for Borel functions.\\

STEP 2: {\sl  It holds that:}
\[
\begin{split}
P(A_1\cap A_2)&=P(A_1)P(A_2)\quad\forall A_1,A_2\in\mathcal
B(\mathbb R).
\end{split}
\]

STEP 3: {\sl Let $h$ be defined by \eqref{eq:h} and for
$A\in\mathcal B(\mathbb R)$ let
\begin{equation}
Z(A,w)=\widetilde{P(A)Wh},
\end{equation}
where $\,\,\widetilde{\mbox{}}\,\,$ denotes the isomorphism
$f\mapsto \widetilde{f}$ from $\mathbf L_2(\mathbb R^n,dx)$ into
$\mathcal
W$ defined by \eqref{eq:iso}. Then $Z$ is a normal field.}\\

Indeed, set $h_0=Wh$. We have
\[
\begin{split}
E\mathbb[ (Z(A_1,\cdot)Z(A_2,\cdot)]&=\langle P(A_1)h_0,
P(A_2)h_0\rangle_{\mathbf L_2(\mathbb R^n,dx)}\\ &=\langle
P(A_1\cap A_2)h_0,h_0\rangle_{\mathbf L_2(\mathbb R^n,dx)}\\
&=\|P(A_1\cap A_2)h_0\|^2_{\mathbf L_2(\mathbb R^n,dx)}\\
&=\|W^*P(A_1\cap A_2)Wh_0\|^2_{\mathbf L_2(\sigma)}\\
&=\|1_{A_1\cap A_2}h_0\|^2_{\mathbf L_2(\sigma)}\\
&=\int_{A_1\cap A_2}\frac{d\sigma(u)}{1+u^2}.
\end{split}
\]

STEP 4: {\sl \eqref{eq:normal} holds.}\\

Indeed,

\[
\begin{split}
X_\sigma(t)&=\widetilde{Q_\sigma(1_{[0,t]})}\\
           &=\widetilde{\sqrt{I+T^2}\widehat{1_{[0,t]}}(T)h_0}\\
           &=\int_{\mathbb
           R}\sqrt{1+u^2}\frac{e^{iut}-1}{u}\widetilde{P(du)h_0}\\
&=\int_{\mathbb
           R}\sqrt{1+u^2}\frac{e^{iut}-1}{u}Z(du,\cdot).
\end{split}
\]

\mbox{}\qed\mbox{}\\



\section{Convolution of measures $\mathcal C$}
\setcounter{equation}{0}
\label{sec:8}
Recall that the family $\mathcal C$ was defined in Definition
\ref{dn:C}. The convolution is not stable in $\mathcal C$. To
verify this, take $n=1$ and $d\sigma=d\lambda$ to be the Lebesgue
measure. Clearly, $d\lambda\in\mathcal C$. On the other hand, we
claim that $d\lambda\star d\lambda\not\in\mathcal C$. Indeed, let
$f\in C_c(\mathbb R)$ (that is, continuous and with support
compact), and such that, moreover
\[
\int_{\mathbb R}f(u)d\lambda(u)=K>0.
\]
We have
\[
\int_{\mathbb R}f(u)d(\lambda\star\lambda)(u)=\iint_{\mathbb
R\times \mathbb R}f(u+v)dudv=K\int_{\mathbb R}dv=\infty.
\]
It follows from Riesz' theorem that $d\lambda\star d\lambda$ is
not
a well defined Borel measure.\\

This example suggests to introduce the class $\mathcal C_b$ which
consists of the positive Borel measures on $\mathbb R^n$ such
that for every $p\in \mathbb N$ there exists $q\in\mathbb N$ and
$C_{pq}>0$ such that
\begin{equation}
\label{eq:rty}
\int_{\mathbb R^n}\frac{d\sigma(u)}{(1+|u+v|^2)^q}\le
\frac{C_{pq}}{(1+|v|^2)^p}.
\end{equation}
\begin{Tm}
It holds that
\[
\mathcal C\star\mathcal C_b\subset\mathcal C.
\]
\end{Tm}

{\bf Proof:} For simplicity we consider the case $n=1$. Let
$\sigma_1\in\mathcal C$ and $\sigma_2\in\mathcal C_b$. There
exists $p\in\mathbb N$ such that
\[
\int_{\mathbb R}\frac{d\sigma_1(u)}{(1+u^2)^p}<\infty
\]
Since $\sigma_2\in\mathcal C_b$, \eqref{eq:rty} is in force for
some $q\in\mathbb N$ and $C_{pq}>0$. Thus:
\[
\begin{split}
\int_{\mathbb R}\frac{d(\sigma_1\star\sigma_2)(w)}{(1+w^2)}=
\iint_{\mathbb
R^2}\frac{d\sigma_1(u)d\sigma_2(v)}{(1+(u+v)^2)^q}\\
\le C_{pq}\int_{\mathbb R}\frac{d\sigma_1(u)}{(1+u^2)^p}<\infty.
\end{split}
\]
\mbox{}\qed\mbox{}\\

\section{An example: The Dirac comb}
\setcounter{equation}{0}
\label{sec:12}
We take
\[
\sigma(u)=\sum_{n\in\mathbb Z}\delta(n-u).
\]
Then, $\mathbf L_2(d\sigma)=\ell_2(\mathbb Z)$. Furthermore:

\begin{Pn}
Let $W$ be an isomorphism between $\ell_2(\mathbb Z)$ onto
$\mathbf L_2(\mathbb R,dx)$, and let $Q$ be defined by
\[
Q\psi=W((\widehat{\psi}(n))_{n\in\mathbb Z}).
\]
Then $Q$ is a bounded operator from $\mathcal S_{\mathbb R}$ into
$\mathbf L_2(\mathbb R,dx)$, and it holds that:
\begin{equation}
\label{clear}
\int_{\mathbb
R}|\widehat{\psi}(u)|^2d\sigma(u)=\int_{\mathbb R}|Q\psi|(x)^2dx
\end{equation}
and $Q^*Q$ is a bounded operator from $\mathcal S_{\mathbb R}$
into $\mathcal S_{\mathbb R}^\prime$ defined by the periodization
operator:
\begin{equation}
(Q^*Q\psi)(x)=\sum_{n\in\mathbb Z}\psi(x+2\pi n).
\label{clear2}
\end{equation}
\end{Pn}

{\bf Proof:} For $\psi\in \mathcal S_{\mathbb R}$, integration by
part shows that the sequence of Fourier coefficients
$(\widehat{\psi}(n))_{n\in\mathbb Z}$ belongs to $\ell_2(\mathbb
Z)$. Therefore, equation \eqref{clear} follows from the definition
of $\sigma$ and from the fact that $W$ is an isomorphism from
$\ell_2(\mathbb Z)$ onto $\mathbf L_2(\mathbb R,dx)$. We now turn
to \eqref{clear2}. Let $\psi,\phi\in\mathcal S_{\mathbb R}$. We
have
\[
\begin{split}
\int_{\mathbb R}(\sum_{n\in\mathbb Z}\psi(x+2\pi
n))\overline{\phi(x)}dx&=\sum_{n\in\mathbb Z} \int_{\mathbb
R}\psi(x+2\pi n)\overline{\phi(x)}dx\\
&=\sum_{n\in\mathbb Z}\frac{1}{2\pi}\int_{\mathbb
R}\widehat{\psi}(u)e^{-2\pi nu i}\overline{\widehat{\phi}(u)}du\\
&=\sum_{n\in\mathbb
Z}\widehat{\psi}(n)\overline{\widehat{\phi}(n)},
\end{split}
\]
where we have used Parseval equality for the second equality, and
Poisson's formula for the third.
\mbox{}\qed\mbox{}\\

\bibliographystyle{plain}
\def\cprime{$'$} \def\lfhook#1{\setbox0=\hbox{#1}{\ooalign{\hidewidth
  \lower1.5ex\hbox{'}\hidewidth\crcr\unhbox0}}} \def\cprime{$'$}
  \def\cprime{$'$} \def\cprime{$'$} \def\cprime{$'$} \def\cprime{$'$}

\end{document}